\documentclass[11pt,reqno]{amsart}
\usepackage{amssymb, latexsym}
\usepackage{graphicx}
\usepackage{fancyhdr}

 \makeatletter
\renewcommand*\subjclass[2][2010]{%
  \def\@subjclass{#2}%
  \@ifundefined{subjclassname@#1}{%
    \ClassWarning{\@classname}{Unknown edition (#1) of Mathematics
      Subject Classification; using '2010'.}%
  }{%
    \@xp\let\@xp\subjclassname\csname subjclassname@#1\endcsname
  }%
}

 \makeatother
\theoremstyle{plain}
\numberwithin{equation}{section}
\newtheorem{thm}{Theorem}[section]
\newtheorem{theorem}[thm]{Theorem}

\newtheorem{corollary}[thm]{Corollary}

\newtheorem{remark}[thm]{Remark}

 \makeatletter
\renewcommand*\subjclass[2][2010]{%
  \def\@subjclass{#2}%
  \@ifundefined{subjclassname@#1}{%
    \ClassWarning{\@classname}{Unknown edition (#1) of Mathematics
      Subject Classification; using '1991'.}%
  }{%
    \@xp\let\@xp\subjclassname\csname subjclassname@#1\endcsname
  }%
}
 \makeatother

\begin{document}
\setcounter{page}{1}

\title[Lucas Type Theorem Modulo Prime 
Powers]{Lucas Type Theorem\\ Modulo Prime Powers} 

\author{Romeo Me\v{s}trovi\'{c}}

\address{Department of Mathematics, Maritime Faculty Kotor,
University of Montenegro\\
Dobrota 36, 85330 Kotor, Montenegro}

\email{romeo@ac.me} 

\begin{abstract} In this note we prove that
 \begin{equation*}
{np^s\choose mp^s+r}\equiv  (-1)^{r-1}r^{-1}(m+1){n\choose m+1}p^s
\pmod{p^{s+1}}
  \end{equation*}
where $p$ is any prime,  $n$, $m$, $s$ and $r$ are nonnegative integers such 
that $n\ge m$, $s\ge 1$, $1\le r\le p^s-1$ and $r$ is not divisible by $p$.
We derive a proof by induction using a multiple application of  Lucas' 
theorem and two basic binomial coefficient identities.
As an application, we prove that a similar congruence for a prime $p\ge 5$
established in 1992 by D. F. Bailey  holds for each prime $p$.
  \end{abstract}
 \maketitle
   

{\renewcommand{\thefootnote}{}

\section{Introduction and Main Result}

In 1878 E.  Lucas \cite{lu} (also see \cite{gr}) proved a remarkable result 
which provides a simple way to compute the binomial coefficient ${n\choose m}$ 
modulo a prime $p$ in terms of the binomial coefficients of the base-$p$ 
digits of $n$ and $m$:  if $n=n_0+n_1p+\cdots +n_sp^s$ and $m=m_0+m_1p+\cdots 
+m_sp^s$ so that $0\le m_i,n_i\le p-1$ for each $i$, then 
   \begin{equation}\label{1.1}
{n\choose m}\equiv \prod_{i=0}^{s}{n_i\choose m_i}\pmod{p}.
   \end{equation}
(with the usual convention that ${0\choose 0}=1$, and 
${l \choose r}=0$ if $l<r$). 
{\em Lucas' theorem} is often formulated in the literature in the 
following equivalent form. If $p$ is a prime, and $a,b,c$ and $d$ are 
nonnegative integers with $a,b\le p-1$, then 
  \begin{equation}\label{1.2}
{cp +a \choose dp+b}\equiv {c \choose d}{a \choose b}
\pmod{p}.
  \end{equation}
In  1990 D. F. Bailey \cite[Theorems 3 and 5]{ba1} proved that 
under the  same assumptions on $a,b,c,d$ for each prime $p\ge 5$
 \begin{equation*}
{cp^f +a \choose dp^f+b}\equiv {c \choose d}{a \choose b}\pmod{p^f}
    \end{equation*}
with $f\in\{2,3\}$. A generalization of  this Lucas-like theorem to every  
prime powers $p^f$ with $p\ge 5$ and $f=2,3,\ldots$ was discovered in 1990 
by K. S. Davis and W. A. Webb \cite{dw1} and independently by A. Granville 
\cite{gr1}.  In 2001 H. Hu and Z.-W. Sun \cite{hs} proved a similar 
congruence  to \eqref{1.2} for generalized binomial coefficients defined in
terms of second order recurrent sequences with initial values 0 and 1.
In  2007 Z.-W. Sun and D. M. Davis \cite{sd} and in 2009 M. Chamberland and 
K. Dilcher \cite{chd} established analogues of Lucas' theorem for certain 
classes of binomial sums. 

Some Lucas' type congruences were  established also by Bailey. 
Namely, in 1991 Bailey \cite[Theorem 4]{ba2} proved by induction on $n\ge 0$ 
that 
  \begin{equation}\label{1.3}  
{np\choose mp+i}\equiv (m+1){n\choose m+1}{p\choose i}\pmod{p^2}
    \end{equation}
where $p$ is a prime, $n$, $m$ and $i$  are nonnegative integers 
with $m\le n$ and $1\le i\le p-1$. 

Applying the congruence \eqref{1.3}, in the same paper \cite[Theorem 5]{ba2} the 
author extended it to the congruence 
 \begin{equation}\label{1.4}
{np^2\choose mp^2+kp+i}\equiv (m+1){n\choose m+1}{p^2\choose kp+i}\pmod{p^3}
    \end{equation}
where $p\ge 5$ is a prime, $n$, $m$, $k$  and $i$  are nonnegative integers
with $m\ge n$, $0\le k\le p-1$ and $1\le i\le p-1$. 

In the next year, proceeding by induction on $s\ge 1$, Bailey 
\cite[Theorem 2.1]{ba3} generalized the congruence \eqref{1.4} modulo higher 
powers of a prime $p\ge 5$. This congruence, extended here for each prime $p$ 
(Corollary 1.2 given below), is obtained as a consequence of the following 
result.

\begin{theorem}\label{1.1}  Let $p$ be any prime,  and let $n$, $m$, $s$ 
and $r$ be nonnegative integers such that $n\le m$, $s\ge 1$, 
$1\le r\le p^s-1$ and $r$ is not divisible by $p$. Then
  \begin{equation}\label{1.5}
{np^s\choose mp^s+r}\equiv  (-1)^{r-1}r^{-1}(m+1){n\choose m+1}p^s
\pmod{p^{s+1}}.
  \end{equation}
$($Here $r^{-1}$ denotes the inverse of $r(\bmod{\,p})$ in the field 
${\bf Z}_p$$)$.
 \end{theorem}

\begin{corollary}\label{1.2} 
{\rm (\cite[Theorem 2.1]{ba3})}.
 Let $p\ge 5$ be a prime and let  $n$, $m$ and $s$  be nonnegative 
integers such that $n\le m$ and $s\ge 1$. Let  $r=\sum_{j=0}^{s-1}a_jp^j$ 
with nonnegative integers $a_j$ such that $1\le a_0\le p-1$ and 
$0\le a_j\le p-1$ for all $j=1,\ldots ,s-1$. Then
  \begin{equation}\label{1.6}
{np^s\choose mp^s+r}\equiv (m+1){n\choose m+1}{p^s\choose r}\pmod{p^{s+1}}.
 \end{equation}
   \end{corollary}

\begin{remark}\label{1.3}
In the proof of Corollary 1.2, using Vandermonde's identity, Bailey  proceed 
by induction on $s$ assuming for the base of induction the cases $s=1$ and 
$s=2$, that is, the congruences {\rm \eqref{1.3}} and {\rm\eqref{1.4}}, 
respectively. Recall that his inductive proof of the congruence 
{\rm\eqref{1.4}} \cite[Theorem 5]{ba1} is based on Vandermonde's identity and
{\it Ljunggren's congruence} {\rm(}see e.g.,  \cite[Theorem 4]{ba1} or 
\cite {gr}{\rm)} which asserts that ${np \choose mp}\equiv {n \choose m}
\,(\bmod{\, p^3})$ for all primes  $p\ge 5$ and nonnegative integers $n$ 
and $m$ with $n\ge m$. Bailey applied the same arguments {\rm(}with 
${np \choose mp}\equiv {n \choose m}\,(\bmod{\, p^2}{\rm)}$ instead of 
Ljunggren's congruence{\rm)} in proof  of the congruence {\rm\eqref{1.3}}
\cite[Theorem 4]{ba1}.
 \end{remark}
In the next section, using only Lucas' theorem and two basic binomial 
coefficient identities, we give an inductive proof of Theorem 1.1.

\section{Proof of Theorem 1.1 and Corollary 1.2}

\begin{proof}[Proof of Theorem $1.1$] First observe that if $n=m$ then since $r\ge 1$,
\eqref{1.5} reduces to the identity $0=0$. Thus, we can assume that $p,n,m$ and $s$
are arbitrary fixed integers satisfying the assumptions of Theorem 1.1 and 
$n\ge m+1\ge 1$. Since by the assumptions, $1\le r\le p^s-1$ and $r$ is not 
divisible by $p$, if $s\ge 2$ we can write $r=kp+i$ with $0\le k\le p^{s-1}-1$ 
and $1\le i\le p-1$, and if $s=1$, then $k=0$ and $r=i$ with $1\le i\le p-1$.

We will prove \eqref{1.5} by induction on $i$ ($=r(\bmod{\,p})$) in the range 
$1\le i\le p-1$. For $i=1$, using the identities 
${a\choose b+1}=\frac{a-b}{b+1}{a\choose b}$ and 
${a\choose b+1}=\frac{a}{b+1}{a-1\choose b}$  with $0\le b\le a-1$, 
we find that 
  \begin{equation}\label{2.1}\begin{split}
&{np^s\choose mp^s+kp+1}=
\frac{(n-m)p^s-kp}{mp^s+kp+1}{np^s\choose mp^s+kp}\\
=&p\cdot\frac{(n-m)p^{s-1}-k}{mp^s+kp+1}
{np^s\choose (mp^{s-1}+k)p}\qquad\qquad\qquad\\
=& p\cdot\frac{(n-m)p^{s-1}-k}{mp^s+kp+1}
\cdot\frac{np^s}{(mp^{s-1}+k)p}{np^s-1\choose (mp^{s-1}+k)p-1}\\
=& p^s\cdot\frac{((n-m)p^{s-1}-k)n}{(mp^s+kp+1)(mp^{s-1}+k)}
{np^s-1\choose (mp^{s-1}+k)p-1}.
  \end{split}\end{equation}
Now we consider two cases.

{\tt Case 1}: $k=0$. Then $r=1$ and for $m=0$ \eqref{1.5} reduces to the identity
$np^s=np^s$. If $m\ge 1$, then the  right hand side of \eqref{2.1} with $k=0$ 
is equal to 
  \begin{equation*}
p^s\cdot\frac{(n-m)n}{(mp^s+1)m}{np^s-1\choose mp^s-1}=
p^s\cdot\frac{(n-m)n}{(mp^s+1)m}
{(np^{s-1}-1)p+(p-1)\choose (mp^{s-1}-1)p+(p-1)}
  \end{equation*}
which  by iterating Lucas' theorem in the form \eqref{1.2} $s$ times and using the 
identity $\frac{(n-m)n}{m}{n-1\choose m-1}=(m+1){n\choose m+1}$ is 
  \begin{eqnarray*}
&\equiv& p^s\cdot\frac{(n-m)n}{m}{np^{s-1}-1\choose mp^{s-1}-1}\pmod{p^{s+1}}\\
&=&p^s\cdot\frac{(n-m)n}{m}{(np^{s-2}-1)p+(p-1)\choose 
(mp^{s-2}-1)p+(p-1)}\pmod{p^{s+1}}\\
&\equiv& p^s\cdot\frac{(n-m)n}{m}{np^{s-2}-1\choose mp^{s-2}-1}\pmod{p^{s+1}}
\equiv\cdots\\
&\equiv& p^s\cdot\frac{(n-m)n}{m}{n-1\choose m-1}=(m+1){n\choose m+1}p^s
\pmod{p^{s+1}}. 
 \end{eqnarray*}
Comparing this with \eqref{2.1} for $k=0$, we find that 
 \begin{equation*}
{np^s\choose mp^s+1}\equiv (m+1){n\choose m+1}p^s\pmod{p^{s+1}}.
  \end{equation*}
This proves \eqref{1.5} with $r=1$ (that is, with $i=1$ and $k=0$).

{\tt Case 2}:  $1\le k\le p^{s-1}-1$.  Then from $r=kp+1\le p^s-1$ we see 
that must be $s\ge 2$. First notice that iterating Lucas' theorem it follows 
immediately that
  \begin{equation}\label{2.2}
{ap^f+c\choose bp^f+d}\equiv {a\choose b}{c\choose d}\pmod{p},
  \end{equation}
where $p$ is a prime, $f,a,b,c$ and $d$ are nonnegative integers such that 
$f\ge 1$, $c\le p^e-1$, $d\le p^e-1$,  and $b\le a$.
In particular, \eqref{2.2} with $c=d=0$ yields
 \begin{equation}\label{2.3}
{ap^e\choose bc^e}\equiv {a\choose b}\pmod{p}.
  \end{equation}
Also notice that for each  prime  $p$ and any integer $j$ such that 
$0\le j\le p-1$  we have
  \begin{equation}\label{2.4}
{p-1\choose j}=\frac{(p-1)(p-2)\cdots (p-j)}{j!}\equiv
\frac{(-1)^jj!}{j!}=(-1)^j\pmod{p}.
   \end{equation}
Take $k=up^l$ where $l\ge 0$ and $u\ge 1$ are nonnegative integers such 
that $u$ is not divisible by $p$. Then since $r=kp+1=up^{l+1}+1\le p^s-1$, 
we see that must be $s\ge 3$, $l\le s-2$ and  $u<p^{s-1-l}$. Taking $k=up^l$ 
for $u=\sum_{j=0}^{s-l-2}u_jp^j$ with $0\le u_j\le p-1$ for each 
$j=0,1,\ldots, s-l-2$ and $u_0\ge 1$ into \eqref{2.1}, using Lucas' theorem, 
\eqref{2.2}, \eqref{2.3} and \eqref{2.4} we find that  
   \begin{equation}\label{2.5}\begin{split}
&{np^s\choose mp^s+up^{l+1}+1}\\
=&p^s\cdot\frac{((n-m)p^{s-1-l}-u)n}{(mp^s+up^{l+1}+1)(mp^{s-l-1}+u)}
{np^s-1\choose (mp^{s-l-1}+u)p^{l+1}-1}\\
\equiv& p^s\cdot\frac{-un}{u}
{(np^{s-l-1}-1)p^{l+1}+(p^{l+1}-1)
\choose (mp^{s-l-1}+u-1)p^{l+1}+(p^{l+1}-1)}\pmod{p^{s+1}}\\
\equiv& -np^s{(np^{s-l-1}-1)p^{l+1}\choose (mp^{s-l-1}+u-1)p^{l+1}}
{p^{l+1}-1\choose p^{l+1}-1}\pmod{p^{s+1}}\\
\equiv& -np^s{(np^{s-l-1}-1)p^{l}\choose (mp^{s-l-1}+u-1)p^{l}}
\pmod{p^{s+1}}\cdots\\
\equiv& -np^s{np^{s-l-1}-1\choose mp^{s-l-1}+u-1}\pmod{p^{s+1}}\cdots\\
=&-np^s{(n-1)p^{s-l-1} +p^{s-l-1}-1\choose mp^{s-l-1}-1+
\sum_{j=0}^{s-l-2}u_jp^j}\\
=&-np^s{(n-1)p^{s-l-1} +\sum_{j=0}^{s-l-2}(p-j)p^j\choose mp^{s-l-1}+(u_0-1)
+\sum_{j=1}^{s-l-2}u_jp^j}\\
\equiv& -np^s{(n-1)p^{s-l-1}\choose mp^{s-l-1}}{p-1\choose u_0-1}
\prod_{j=1}^{s-l-2}{p-1\choose u_j}\pmod{p^{s+1}}\equiv\cdots\\
\equiv& -np^s{n-1\choose m}(-1)^{-1+\sum_{j=0}^{s-l-2}u_j}\pmod{p^{s+1}}\\
\equiv& n{n-1\choose m}(-1)^{u}p^s\pmod{p^{s+1}}\\
\equiv& (m+1){n\choose m+1}(-1)^{r-1}p^s\pmod{p^{s+1}}
   \end{split}\end{equation}
(the last two congruences are clearly satisfied since for odd prime $p$, 
$\sum_{j=0}^{s-l-2}u_j\equiv u(\bmod{\,2})$, and hence 
$r-1=up^{l+1}\equiv u(\bmod{\,2})$, while for $p=2$ we have 
$(-1)^t\equiv 1(\bmod{\,2})$ for each integer $t$).
The congruence \eqref{2.5} coincides with \eqref{1.5} for $r=up^{l+1}+1$. 
This concludes the proof of the base of induction ($i=1$).

Now suppose that the congruence \eqref{1.5} holds for each $r=kp+i$ with 
$0\le k\le p^{s-1}-1$ and some fixed $i$ with $1\le i\le p-2$; that is 
  \begin{equation}\label{2.6}
{np^s\choose mp^s+kp+i}\equiv (-1)^{kp+i-1}(kp+i)^{-1}(m+1){n\choose m+1}
p^s\pmod{p^{s+1}}.
  \end{equation}
Then using the identity ${a\choose b+1}=\frac{a-b}{b+1}{a\choose b}$  
 with $0\le b\le a$ and \eqref{2.6}, we find that 
  \begin{eqnarray*}
&&{np^s\choose mp^s+kp+i+1}=
\frac{(n-m)p^s-kp-i}{mp^s+kp+i+1}{np^s\choose mp^s+kp+i}\\
&\equiv& \frac{(n-m)p^s-kp-i}{mp^s+kp+i+1}
(-1)^{kp+i-1}(kp+i)^{-1}(m+1){n\choose m+1}p^s\pmod{p^{s+1}}\\
&\equiv& \frac{-i}{kp+i+1}
(-1)^{kp+i-1}i^{-1}(m+1){n\choose m+1}p^s\pmod{p^{s+1}}\\
&=& (-1)^{kp+i}(kp+i+1)^{-1}(m+1){n\choose m+1}p^s\pmod{p^{s+1}}.
  \end{eqnarray*}
This proves \eqref{1.5} with $r\equiv i+1(\bmod{\, p})$, which completes 
proof of Theorem 1.1.
\end{proof}

\begin{proof}[Proof of Corollary $1.2$] Taking $n=1$ and $m=0$ into the congruence 
\eqref{1.5} of Theorem 1.1, for all $r$ such that $1\le r\le p^s-1$ and $r$ is not 
divisible by $p$ we get
  \begin{equation*}
{p^s\choose r}\equiv  (-1)^{r-1}r^{-1}p^s\pmod{p^{s+1}}.
  \end{equation*}
Comparing this with \eqref{1.5}, we immediately obtain \eqref{1.6}.
\end{proof}


\begin{thebibliography}{99} 
\bibitem{ba1} D. F. Bailey,  \emph{Two $p^3$ variations of Lucas' theorem}, 
J. Number Theory  35 (1990), 208--215. 

\bibitem{ba2} D. F. Bailey,  \emph{Some binomial coefficient congruences, 
 Appl. Math. Lett.} \textbf{4}, No. 4 (1991), 1--5. 

\bibitem{ba3} D. F. Bailey,  \emph{More binomial coefficent congruences}, 
 Fibonacci Quart.  \textbf{30}, No. 2 (1992), 121--125. 

\bibitem{dw1} K. S. Davis and W. A. Webb, \emph{Lucas' theorem 
for  prime powers}, European J. Combin.   \textbf{11} (1990), 229--233. 

\bibitem{chd} M. Chamberland and K. Dilcher, \emph{A binomial sum related 
to Wolstenholme's theorem}, J. Number Theory  \textbf{129} (2009), 
2659--2672.

\bibitem{gr} A. Granville,  \emph{Arithmetic properties of binomial 
coefficients. I. Binomial coefficients modulo prime powers}, in Organic mathematics 
(Burnaby, BC, 1995), CMS Conf. Proc., 20, American Mathematical Society, 
Providence, RI, 1997, 253--276. 

\bibitem{gr1} A. Granville, \emph{Zaphod Beeblebrox's brain and the 
fifty-ninth row of Pascal's triangle}, Amer. Math. Monthly  \textbf{99} 
(1992), 318--331.  


\bibitem{hs} H. Hu and Z.-W. Sun, \emph{An extension of Lucas' theorem},
 Proc. Amer. Math. Soc.  \textbf{129} (2001), 3471--3478.

\bibitem{lu} \'{E}. Lucas, \emph{Sur les congruences des nombres 
eul\'{e}riens et des coefficients diff\'{e}rentiels des fonctions 
trigonom\'{e}triques, suivant un module premier}, 
 Bull. Soc. Math. France  \textbf{6} (1877--1878), 49--54.

\bibitem{sd} Z.-W. Sun and D. M. Davis, \emph{Combinatorial congruences 
modulo prime powers},  Trans. Amer. Math. Soc.  \textbf{359} (2007), 5525--5553. 
 \end{thebibliography}
\end{document}